\documentclass[12pt,leqno]{article}

\usepackage{amsmath}              
\usepackage[mathscr]{eucal}       
\usepackage{amssymb}              
\usepackage{amsbsy}               
\usepackage{amsthm}               
\usepackage{fancyhdr}             
\usepackage{color}		    

\renewcommand{\headrule}{\vbox to 0pt{\hbox to\headwidth{%
	\textcolor[cmyk]{0,0,.8,.2}{\rule{\headwidth}{1pt}}}\vss}}

\def\hours{\begingroup
	\newcount\hour
	\hour=\time
	\divide\hour by 60 
	\ifnum\hour < 1 12\fi
	\ifnum\hour > 0
		\the\hour
	\fi
	\endgroup}
\def\minutes{\begingroup
	\newcount\min
	\min=\time
	\divide\time by 60
	\multiply\time by 60
	\advance\min by -\time 
		\ifnum\min < 10 0\fi \the\min 
	\endgroup }
\newcommand{\timenow}{\hours:\minutes}

\def\APtimenow{\begingroup
	\ifnum\time < 60 12:\minutes\ A.M
	\else 
		\ifnum\time > 720 \advance\time by -720 \timenow\ P.M.
		\else \timenow\ A.M.
		\fi
	\fi
	\endgroup}

\newcommand{\theAuthor}{}
\newcommand{\theTitle}{}
\newcommand{\dayFinished}{}

\newcommand{\niceheader}[3]{%
	\setAuthor{#1}
	\setTitle{#2}
	\setDayFinished{#3}
	\vspace{-1cm}\noindent  {\color[cmyk]{0,0,.8,.2}{ \rule[10pt]{\textwidth}{1.4pt} }}
	\centerline{ \large\sf{\theTitle} }
	\centerline{\it{\sf\theAuthor}}
	\centerline{\it{\sf\dayFinished}}
	\noindent  {\color[cmyk]{0,0,.8,.2}{ \rule[-2pt]{\textwidth}{1.4pt} }}
	\thispagestyle{plain} 
	}

\newcommand{\setAuthor}[1]{\def\theAuthor{#1}}
\newcommand{\setTitle}[1]{\def\theTitle{#1}}
\newcommand{\setDayFinished}[1]{\def\dayFinished{#1}}

\newcommand{\genericFooter}{%
	\label{doc:end}
	\noindent\rule{\textwidth}{1.4pt}
	\rightline{\tiny{Last compiled at \APtimenow\ on \today }}
}

\AtEndDocument{\genericFooter}

\definecolor{dw04blue}{cmyk}{1,.5,0,.6}
\definecolor{dw04purple}{cmyk}{.3,.9,0,.3}
\definecolor{dw04red}{cmyk}{0,.8,.7,.5}
\definecolor{dw04orange}{cmyk}{0,.4,1,.1}
\definecolor{dw04yellow}{cmyk}{0,0,1,.1}
\definecolor{dw04green}{cmyk}{.9,0,.9,.6}
\definecolor{dw04black}{cmyk}{0,0,0,1}

\newtheoremstyle{dw-defn}
  {}
  {}
  {}
  {}
  {\sc}
  {:}
  {.4em}
  {}

\newtheoremstyle{dw-plain}
  {}
  {}
  {}
  {}
  {\sc}
  {:}
  {.4em}
  {}

\theoremstyle{dw-plain}
\swapnumbers
\newtheorem{thm}{\textcolor{dw04purple}{Theorem}}[section]
\newtheorem{lem}[thm]{\textcolor{dw04purple}{Lemma}}
\newtheorem{prop}[thm]{\textcolor{dw04purple}{Proposition}}
\newtheorem{cor}[thm]{\textcolor{dw04purple}{Corollary}}

\theoremstyle{dw-defn}

\newcommand{\union}{\cup}

\newcommand{\naturals}{\mathbb{N}}

\renewcommand{\subset}{\subseteq}

\newcommand{\focus}[1]{\textcolor[cmyk]{.2,.8,1,.2}{{\bf{#1}}}}
\newcommand{\stend}{\hfill $\lozenge$}

\newcommand{\stigma}{\sigma^*}
\newcommand{\e}{\varepsilon}
\newcommand{\floor}[1]{\left\lfloor#1\right\rfloor}
\newcommand{\ceil}[1]{\left\lceil#1\right\rceil}

\newcommand{\separator}{\rule[6pt]{50pt}{.6pt}}

\begin{document}
\niceheader{Dan Warren}{Optimal Packing Behavior of some $2$-block Patterns}{March 18, 2004}
\begin{abstract}
In this paper, a result of Albert, Atkinson, Handley, Holton, and
Stromquist (Proposition 2.4 of \cite{stromquist-2001}) which characterizes the optimal
packing behavior of the pattern 1243 is generalized in two directions.
The packing densities of layered patterns of type $(1^\alpha,\alpha)$
and $(1,1,\beta)$ are computed.
\end{abstract}
\centerline{\small{{\bf{AMS Subject Classifications:}} 05A15, 05A16}}
\centerline{\small{{\bf{Keywords:}} pattern containment, permutations, layered permutations, packing density.}}

\section{Definitions and Notation} 

\begin{figure}[b]
\begin{center}
\begin{picture}(150,150)(0,0) 

  \definecolor{dw04blue}{cmyk}{1,.5,0,.1}
  \definecolor{dw04red}{cmyk}{0,.8,.7,.2}
  \definecolor{dw04black}{cmyk}{0,0,0,1}

	\color{dw04black}
	\put(0,150){\line(1,0){150}}
	\put(0,150){\line(0,-1){150}}
	\put(0,0){\line(1,0){150}}
	\put(150,150){\line(0,-1){150}}

	\color{dw04blue}
	\put(53,2){\line(0,1){147}}

	\put(87,2){\line(0,1){146}}

	\put(2,52){\line(1,0){146}}

	\put(2,83){\line(1,0){145}}

	\color{dw04red}
	\put(12,39){\circle*{6}}

	\put(26,27){\circle*{6}}

	\put(38,15){\circle*{6}}

	\put(61,75){\circle*{6}}

	\put(74,64){\circle*{6}}

	\put(94,139){\circle*{6}}

	\put(108,126){\circle*{6}}

	\put(123,111){\circle*{6}}

	\put(137,97){\circle*{6}}

\end{picture} 
\end{center}
\caption{The permutation 321549876 is layered, with layer sizes $(3,2,4)$ }
\label{fig:layPermEx}
\end{figure}
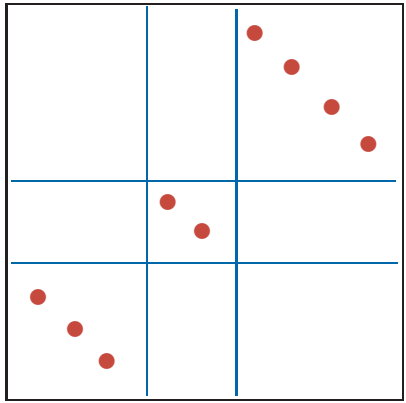

Let $\sigma$ be a permutation on $[1..n]$, and let $\tau$ be a pattern on $[1..m]$,
for some $m \leq n$.  We say that an $m$-subset $S \subset [1..n]$ is an
\focus{occurrence} of $\tau$ in $\sigma$ if the restriction of $\sigma$ to $S$ is
isomorphic to $\tau$, that is, they are in the same linear order.  
There have been
two major areas of study within the realm of pattern containment, avoidance and packing.  
The usual aim of the pattern avoidance problem is to enumerate or otherwise characterize the permutations
on $[1..n]$ which avoid a specific pattern or set of patterns.  In contrast, the study 
of packing patterns into permutations, which was begun in the early 1990s as an offshoot of 
the study of pattern avoidance, aims at the opposite question:  Given a fixed pattern $\tau$
(or set of patterns, as in \cite{stromquist-2001}), how must we structure a permutation
$\sigma$ of length $n$ so that it has the maximal number of occurrences of $\tau$?
Much of the founding work on the subject was done by Alkes Price, in his Ph.D. dissertation (\cite{price}),
and most of the literature on the subject addresses the class of \focus{layered}
permutations, which we will define below.  As \cite{price} is not widely available,
the reader is referred to \cite{bona-pbook}, which addresses a broad range of 
topics related to pattern containment.

Adopting the notation of Price, we define a \focus{layer} in $\sigma$ to be a contiguous
decreasing subsequence of consecutive integers.  A pattern is called \focus{layered}
if it consists of an increasing sequence of (disjoint) layers (a layered permutation is
shown in figure \ref{fig:layPermEx}).  

It has long been thought that the optimal packing behavior of patterns
containing a contiguous block of $r$ layers of size 1 was analagous to
that of patterns having one layer of size $r$ in the corresponding
place.
In \cite{stromquist-2001}, the authors were the first to address nontrivial
patterns having multiple layers of size 1 in a row, with the
example of the pattern 1243.  The authors of that paper showed
that the optimal packing behavior of 1243 is analogous (in its simplicity) to
that of 2143, that is, there is a single increasing block,
followed by a single decreasing block.  In this paper, we
provide two slight generalizations of their work, to patterns
having layer sizes $(1^\alpha,\alpha)$ and patterns
having layer sizes $(1,1,\beta)$, both of which display similar behavior.  

To deal efficiently with contiguous blocks of layers of size 1, we
define an \focus{antilayer} to be a contiguous {\em{increasing}} subsequence of
consecutive integers, that is, a contiguous sequence of layers of size 1.  The more general 
term \focus{block} will be understood
to mean either a layer or an antilayer.  We will use the term \focus{isolated point}
to describe a layer of size 1 between two layers of larger sizes.
Since \cite{price} settles the question of optimal packing behavior for patterns having 2 layers,
the next logical step is to attack the slightly more general class of
patterns having two blocks.  The only such patterns not characterized by
Price's work are those which consist of one antilayer and one
layer [each of size at least 2].

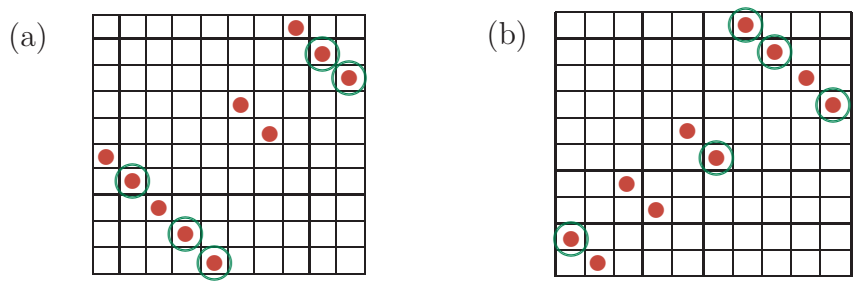
\begin{figure}[t]
\begin{center}

\begin{picture}(350,120)(20,20) 

  \definecolor{dw04red}{cmyk}{0,.8,.7,.2}
  \definecolor{dw04green}{cmyk}{.9,0,.8,.3}
  \definecolor{dw04black}{cmyk}{0,0,0,1}

	\color{dw04black}
	\put(49,111){\line(1,0){103}}
	\put(49,102){\line(1,0){103}}
	\put(49,92){\line(1,0){103}}
	\put(49,82){\line(1,0){103}}
	\put(49,72){\line(1,0){103}}
	\put(49,62){\line(1,0){103}}
	\put(49,53){\line(1,0){103}}
	\put(49,43){\line(1,0){103}}
	\put(49,33){\line(1,0){103}}
	\put(49,23){\line(1,0){103}}
	\put(49,13){\line(1,0){103}}
	\put(49,111){\line(0,-1){98}}
	\put(59,111){\line(0,-1){98}}
	\put(69,111){\line(0,-1){98}}
	\put(79,111){\line(0,-1){98}}
	\put(90,111){\line(0,-1){98}}
	\put(100,111){\line(0,-1){98}}
	\put(110,111){\line(0,-1){98}}
	\put(121,111){\line(0,-1){98}}
	\put(131,111){\line(0,-1){98}}
	\put(141,111){\line(0,-1){98}}
	\put(152,111){\line(0,-1){98}}

	\color{dw04black}
	\put(224,112){\line(1,0){112}}
	\put(224,102){\line(1,0){112}}
	\put(224,92){\line(1,0){112}}
	\put(224,82){\line(1,0){112}}
	\put(224,72){\line(1,0){112}}
	\put(224,62){\line(1,0){112}}
	\put(224,52){\line(1,0){112}}
	\put(224,42){\line(1,0){112}}
	\put(224,32){\line(1,0){112}}
	\put(224,22){\line(1,0){112}}
	\put(224,12){\line(1,0){112}}
	\put(224,112){\line(0,-1){100}}
	\put(235,112){\line(0,-1){100}}
	\put(246,112){\line(0,-1){100}}
	\put(257,112){\line(0,-1){100}}
	\put(268,112){\line(0,-1){100}}
	\put(280,112){\line(0,-1){100}}
	\put(291,112){\line(0,-1){100}}
	\put(302,112){\line(0,-1){100}}
	\put(313,112){\line(0,-1){100}}
	\put(324,112){\line(0,-1){100}}
	\put(336,112){\line(0,-1){100}}

	\color{dw04black}
	\put(17,99){(a)}

	\color{dw04black}
	\put(198,100){(b)}

	\color{dw04red}
	\put(54,57){\circle*{6}}

	\color{dw04red}
	\put(64,48){\circle*{6}}

	\color{dw04red}
	\put(74,38){\circle*{6}}

	\color{dw04red}
	\put(84,28){\circle*{6}}

	\color{dw04red}
	\put(95,17){\circle*{6}}

	\color{dw04red}
	\put(105,77){\circle*{6}}

	\color{dw04red}
	\put(116,66){\circle*{6}}

	\color{dw04red}
	\put(126,106){\circle*{6}}

	\color{dw04red}
	\put(136,96){\circle*{6}}

	\color{dw04red}
	\put(146,87){\circle*{6}}

	\color{dw04green}
	\put(64,48){\circle{12}}
	\put(64,48){\circle{13}}

	\put(84,28){\circle{12}}
	\put(84,28){\circle{13}}

	\put(95,17){\circle{12}}
	\put(95,17){\circle{13}}

	\put(136,96){\circle{12}}
	\put(136,96){\circle{13}}

	\put(146,87){\circle{12}}
	\put(146,87){\circle{13}}

	\color{dw04red}
	\put(230,26){\circle*{6}}

	\put(240,17){\circle*{6}}

	\put(251,47){\circle*{6}}

	\put(262,37){\circle*{6}}

	\put(274,67){\circle*{6}}

	\put(285,57){\circle*{6}}

	\put(296,107){\circle*{6}}

	\put(307,97){\circle*{6}}

	\put(319,87){\circle*{6}}

	\put(329,77){\circle*{6}}

	\color{dw04green}
	\put(230,26){\circle{12}}
	\put(230,26){\circle{13}}

	\put(285,57){\circle{12}}
	\put(285,57){\circle{13}}

	\put(296,107){\circle{12}}
	\put(296,107){\circle{13}}

	\put(307,97){\circle{12}}
	\put(307,97){\circle{13}}

	\put(329,77){\circle{12}}
	\put(329,77){\circle{13}}

\end{picture} 

\end{center}
\caption{(a) Notice that in each occurrence of the pattern $\tau=32154$, each
layer of $\tau$ must be contained in a layer of the larger permutation $[5,4,3,2,1,7,6,10,9,8]$.  An occurrence
of $\tau$ is circled.  (b) Observe
how the antilayer in the pattern $\tau=12543$ may `climb' layers in the larger permutation
$[2,1,4,3,6,5,10,9,8,7]$ rather than being contained in a single antilayer.
An occurrence of $\tau$ is circled. }
\separator\vspace{-.5cm}
\label{fig:climbing}
\end{figure}

While the concepts of a layer and an antilayer look virtually the same, 
they unfortunately must
be handled in quite different ways, given the traditional approach to the
problem.  Due to the work of Price and Stromquist, when $\tau$ is a layered pattern, 
we are able to restrict our search for $\tau$-maximal permutations of $[1..n]$ to 
the much smaller class of layered permutations of $[1..n]$.  Note, however, 
that in an occurrence of the pattern $\tau$ in a permutation $\sigma$, any layer of $\tau$
must be contained in a layer of $\sigma$, but antilayers in $\tau$ need not
be contained in antilayers of $\sigma$, as they can `climb' a list of several
layers, having one element in each (see Figure \ref{fig:climbing}).  This simple problem is enough to
make the computation notably more difficult, so that a slightly more
delicate argument is required to prove a conjecture which at the outset
seems obvious.

For convenience of notation, let $\tau_{\alpha,\beta}$ be the layered
pattern having layers $(1^\alpha,\beta)$.
We will adopt some of Price's notation for the counting of occurrences: 
for $\sigma \in S_n$ and any $\tau \in S_m$, $m<n$, let 
$g(\tau,\sigma)$
denote the number
of occurrences of $\tau$ in $\sigma$.  Let 
	$$g(\tau,n) = \max_{\sigma \in S_n} g(\tau,\sigma),$$
and say that $\sigma \in S_n$ is $\tau$-maximal if $g(\tau,\sigma) = g(\tau,n)$.

A theorem attributed to Galvin (reproduced as Theorem 2.1 of \cite{price}) states that
the sequence $\left( \frac{g(\tau,n)}{\binom{n}{m}} \right)_{n \geq m}$ is in fact 
decreasing in $n$, so that the limit
	$$\lim_n \frac{g(\tau,n)}{\binom{n}{m}}$$
exists.  This limit, which we will denote $\rho(\tau)$, is called the \focus{packing density} of
$\tau$.  We will compute the packing density of each class of patterns we explore.

\section{The Layered Pattern $\tau_{\alpha,\alpha}$}

The idea in this section is that given the packing behavior of a small pattern, we may
determine the packing densities for larger, similar patterns
by exploiting the similarity in structure.

\begin{prop} The structure of the maximizing permutation of size $2n$ for the 
pattern 
	$$\tau_{\alpha,\alpha}:=[1,2,\dots \alpha,2\alpha,2\alpha-1,\dots \alpha+1]$$
 is invariant of $\alpha$, that is, the $\tau_{\alpha,\alpha}$-maximizing
pattern of length $2n$ is of the form $[1,2,\dots n,2n,2n-1,\dots, n+1]$ for all $\alpha$.
\begin{proof}
The case $\alpha=2$ was proven in Proposition 2.4 of \cite{stromquist-2001}; it follows that the maximal number of
$1243$s in a pattern of length $2n$ is $\binom{n}{2}^2$.  In each step, we will show that
if the pattern $[1,2,\dots,n,2n,2n-1,\dots,n+1]$ is not $\tau_{\alpha,\alpha}$-maximal, then in fact it cannot
be $\tau_{2,2}$-maximal, which would contradict the known result for $1243$.
We will first prove the $\alpha=3$ case.
Consider $\tau_{3,3} = 123654$, and suppose that $g(\tau_{3,3},2n) \geq \binom{n}{3}^2 + 1$.  Let $\sigma \in S_{2n}$
be $\tau_{3,3}$-maximal.  In each $\tau_{3,3}$, there are $9=\binom{3}{2}^2$ instances of $1243$s.  Suppose that
$\sigma$ is an instance of $1243$ in $\sigma$.  Now, each of the remaining $2n-4$ elements, if in a 
$\tau_{3,3}$ containing $\sigma$, can be in either side, the $123$, or the $654$, but not both, so 
there are $\ell \cdot (2n-4-\ell)$ ways to form a $\tau_{3,3}$ for some $\ell$.  Since the largest value of
the expression $\ell\cdot(2n-4-\ell)$ occurs when $\ell=n-2$, it follows that the number of $\tau_{3,3}$s containing
a $1243$ is at most $(n-2)^2$.  Hence, if $\sigma$ has at least $\binom{n}{3}^2 + 1$
occurrences of $\tau_{3,3}$, then it must have at least
	$$\frac{\binom{3}{2}^2}{(n-2)^2}\left[ \binom{n}{3}^2 + 1 \right] > \binom{n}{2}^2$$
occurrences of $1243$, contradicting \cite{stromquist-2001}.

The general case is quite similar:
 Suppose
that $g(\tau_{\alpha,\alpha},2n) \geq \binom{n}{k}^2 + 1$ and let $\sigma \in S_{2n}$ be a permutation having $g(\tau_{\alpha,\alpha},2n)$
occurrences of $\tau_{\alpha,\alpha}$.  We will again count the number of $1243$s in $\sigma$:  In each
occurrence of $\tau_{\alpha,\alpha}$, there are $\binom{\alpha}{2}^2$ occurrences of $1243$.  Given a particular
occurrence $\tau_0$ of $1243$, we need to determine the maximum number of $\tau_{\alpha,\alpha}$s which could contain $\tau_0$;
however, for each $\tau_{\alpha,\alpha}$ containing $\tau_0$, the other $2\alpha-4$ of its elements must be
a $\tau_{\alpha-2,\alpha-2}$ in the remaining $2n-4$ elements of $\sigma$.  That is, the number of occurrences of $\tau_{\alpha,\alpha}$ which
can contain $\tau_0$ is bounded above by $g(\tau_{\alpha-2,\alpha-2},2n-4)$.  By induction, there are at most
$\binom{n-2}{\alpha-2}^2$ of these. 
Hence, the number of $1243$s in $\sigma$ is at least
	$$\frac{\binom{\alpha}{2}^2}{\binom{n-2}{\alpha-2}^2}\left[ \binom{n}{\alpha}^2 + 1 \right] = \binom{n}{2}^2 + \e$$
where $\e > 0$, which again contradicts the result of \cite{stromquist-2001}.  It follows that $g(\tau_{\alpha,\alpha},2n) \leq \binom{n}{\alpha}^2$
for each $n$, which means that the structure of a $\tau_{\alpha,\alpha}$-maximizing permutation on $[1..2n]$ is
of the form $[1,\dots,n,2n,2n-1,\dots n+1]$.
\end{proof}
\end{prop}

\begin{cor}  The packing density of the pattern $\tau_{\alpha,\alpha}$
above is
	$$\rho(\tau_{\alpha,\alpha}) = \frac{\binom{2\alpha}{\alpha}}{2^{2\alpha}}.$$
\begin{proof}
The result follows from taking the limit in $n$ of the fraction
	$$g(\tau_{\alpha,\alpha},2n) = \frac{\binom{n}{\alpha}\binom{n}{\alpha}}{\binom{2n}{2\alpha}}.$$
\end{proof}
\end{cor}

These results are nice enough, but they still only apply to a very narrow class of patterns.  The next step forward from
here is to try to prove a statement about patterns of the same format, but with 2 blocks 
being permitted to have different sizes.

\section{The Layered Pattern $\tau_{2,\beta}$}

\noindent As the optimal packing behavior of the layered pattern $\tau_{\alpha,\beta}$ 
does not in general adhere to the same degree of symmetry as the case $\alpha=\beta$, the ideas of the previous section do not so easily
lend themselves to the general case, so for the duration of this section 
we must cease to rely on known results and
do some computation of our own.

We will later need the following technical lemma:
\begin{lem}\label{lem:tech}
Let $k,\ell,m,n \in \naturals$ s.t. $k < \ell \leq m \leq n$.
Then, we have
	$$\binom{n}{k} \binom{m}{\ell} \leq \binom{n}{\ell} \binom{m}{k}.$$
That is, we can have more combinations if bigger sets choose the bigger subsets.
\begin{proof}
Although I am sure there is an elegant combinatorial proof, the following straight computation will do,
as this result is only developed here as a tool for later use.  We have
	$$\frac{ \binom{m}{\ell} }{ \binom{m}{k} } = \frac{k!}{\ell!} \cdot (m-l)(\dots )(m-\ell+1)$$
and
	$$\frac{ \binom{n}{\ell} }{ \binom{n}{k} } = \frac{k!}{\ell!} \cdot (n-l)(\dots)(n-\ell+1)$$
so that the result follows from cross-multiplying.
\end{proof}
\end{lem}

For the remainder of this section, let $\beta \geq 3$ and let $\sigma \in S_n$ be a permutation which is $\tau_{2,\beta}$-maximal.  

An early result due to Stromquist provides us with the first important assumption we can make about $\sigma$,
namely that it is layered.

\begin{thm}[Theorem 2.2 of \cite{price}]
Let $\tau \in S_m$ be a layered pattern.  Then, for each $n \geq m$, we have
	$$g(\tau,n) = \max \left\{ g(\tau,\sigma) :: \sigma \in S_n \right\} = \max \left\{ g(\tau,\sigma) :: \sigma \in S_n \text{ is layered} \right\}.$$
\stend
\end{thm}


\begin{lem}\label{lem:isolated}
We can 
assume without loss of generality that $\sigma$ has no isolated points, i.e. layers
of size $1$ between two layers of size $\geq 2$.
\begin{proof}
Suppose we have an isolated point
in between layers $L_i$ and $L_{i+1}$.  Switching the positions of $L_i$ and the isolated 
point does not decrease the number of occurrences of $\tau$, so we may as well do that.
Since we can always move an isolated point to the left of a layer, eventually all isolated
points will be soaked up by antilayers.  We can assume, then, that $\sigma$ has no
isolated points.
\end{proof}
\end{lem}

It follows now that $\sigma$ is a sequence of layers and antilayers, with no isolated points.
In fact, a more general statement is true by the same logic:

\begin{figure}[b]
\begin{center}
\begin{picture}(150,150)(20,20) 

  \definecolor{dw04blue}{cmyk}{1,.5,0,.1}
  \definecolor{dw04purple}{cmyk}{.3,.9,0,.2}
  \definecolor{dw04orange}{cmyk}{0,.4,1,0}
  \definecolor{dw04black}{cmyk}{0,0,0,1}

	\color{dw04black}
	\put(0,150){\line(1,0){150}}
	\put(0,150){\line(0,-1){150}}
	\put(0,0){\line(1,0){150}}
	\put(150,150){\line(0,-1){150}}

	\color{dw04blue}
	\put(28,2){\line(0,1){146}}

	\color{dw04blue}
	\put(54,2){\line(0,1){148}}

	\color{dw04blue}
	\put(121,1){\line(0,1){148}}

	\color{dw04blue}
	\put(1,27){\line(1,0){148}}

	\color{dw04blue}
	\put(2,58){\line(1,0){148}}

	\color{dw04blue}
	\put(2,123){\line(1,0){148}}

	\color{dw04orange}
	\put(1,2){\line(1,1){25}}

	\color{dw04orange}
	\put(28,54){\line(1,-1){25}}

	\color{dw04orange}
	\put(120,148){\line(6,-5){27}}

	\color{dw04purple}
	\put(11,9){$A_1$}

	\color{dw04purple}
	\put(35,37){$L_1$}

	\color{dw04purple}
	\put(126,132){$L_k$}

	\color{dw04purple}
	\put(75,75){.}

	\color{dw04purple}
	\put(90,90){.}

	\color{dw04purple}
	\put(107,105){.}

\end{picture} 

\end{center}
\caption{The form of $\sigma$, as shown in Lemma \ref{lem:layerform} }
\label{fig:layerform}
\end{figure}
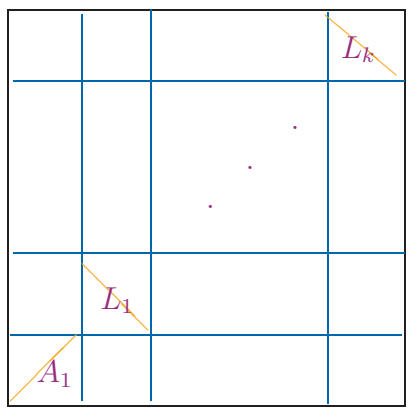

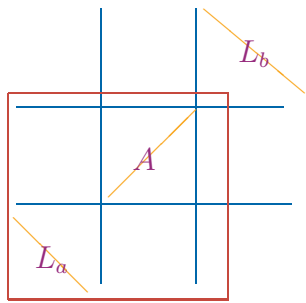
\begin{figure}[t]
\begin{center}

\begin{picture}(120,120)(20,20) 

  \definecolor{dw04blue}{cmyk}{1,.5,0,.1}
  \definecolor{dw04purple}{cmyk}{.3,.9,0,.2}
  \definecolor{dw04red}{cmyk}{0,.8,.7,.2}
  \definecolor{dw04orange}{cmyk}{0,.4,1,0}
  \definecolor{dw04black}{cmyk}{0,0,0,1}

	\color{dw04blue}
	\put(37,113){\line(0,-1){104}}

	\color{dw04blue}
	\put(73,113){\line(0,-1){104}}

	\color{dw04blue}
	\put(5,39){\line(1,0){104}}

	\color{dw04blue}
	\put(5,76){\line(1,0){101}}

	\color{dw04orange}
	\put(4,34){\line(1,-1){28}}

	\color{dw04orange}
	\put(40,42){\line(1,1){33}}

	\color{dw04orange}
	\put(76,113){\line(6,-5){38}}

	\color{dw04purple}
	\put(12,15){$L_a$}

	\color{dw04purple}
	\put(89,93){$L_b$}

	\color{dw04purple}
	\put(49,52){$A$}

	\color{dw04red}
	\put(2,81){\line(1,0){83}}
	\put(2,81){\line(0,-1){78}}
	\put(2,3){\line(1,0){83}}
	\put(85,81){\line(0,-1){78}}

\end{picture} 

\end{center}
\caption{There can be no occurrences of $\tau_{2,\beta}$ inside the boxed area, so it
will not destroy any occurrences of $\tau_{2,\beta}$ to switch the positions of $L_a$
and $A$.}
\label{fig:laswap}
\separator
\end{figure}

\begin{lem}\label{lem:layerform}
We may assume that $\sigma$ consists of a single antilayer $A_1$ followed by a
list of (nontrivial) layers $L_1,\dots,L_k$, as shown in Figure \ref{fig:layerform}.
\begin{proof}
Suppose that in the middle of the permutation $\sigma$ we have an antilayer
between two layers.  Similar to the proof of Lemma \ref{lem:isolated}, if we switch the 
positions of the antilayer and the layer to its left (see figure \ref{fig:laswap}), all the occurrences of $\tau_{2,\beta}$
that were originally there are left intact, and we create more as long as the antilayer
has length at least $2$ and the left layer has length at least $\beta$.  Thus,
in a permutation $\sigma$, we may move all antilayers to the left of all layers without
losing any occurrences of $\tau_{2,\beta}$,
achieving the desired layer/antilayer pattern.
\end{proof}
\end{lem}

\begin{lem}\label{lem:sorting}
The layers $L_1,\dots,L_k$ may be assumed to be in nondecreasing order by size.
\begin{proof}
Suppose $L_i$ and $L_{i+1}$ are adjacent layers, and that $L_i$ is larger than $L_{i+1}$.
What happens if we switch the positions of the two layers?  Unless an occurrence of $\tau_{2,\beta}$
has its $\beta$-layer in one of $L_i$ or $L_{i+1}$, it will be preserved, just moved, so
we need only worry about those occurrences of $\tau_{2,\beta}$ whose
$\beta$-layer is contained in either $L_i$ or $L_{i+1}$.
Let $y = |A_1| + \sum_{j=1}^{i-1}|L_i|$ be the number of elments of $\sigma$ to the
left of $L_i$, and let $\alpha$ be the number of increasing $2$-sequences
in this range.  Then, the number of occurrences of $\tau_{2,\beta}$ which have
a $\beta$-layer in $L_i$ or $L_{i+1}$ is
	$$ \alpha \binom{|L_i|}{\beta} + \alpha \binom{|L_{i+1}|}{\beta}
		+ y \cdot |L_i| \binom{|L_{i+1}|}{\beta}.$$
If we swap the layer-lengths, the number of these occurrences becomes
	$$ \alpha \binom{|L_{i+1}|}{\beta} + \alpha \binom{|L_{i}|}{\beta}
		+ y \cdot |L_{i+1}| \binom{|L_{i}|}{\beta}.$$
Of course the sum of the first two terms remains the same.  However, since $L_{i} > L_{i+1}$, 
the final term is certainly larger after the switch, by our technical lemma.  It follows that if $|L_i| > |L_{i+1}|$,
the permutation $\sigma$ cannot be $\tau_{2,\beta}$-maximal.  Since this statement
holds for all $i \in [1..k-1]$, we can assume that if $\sigma$ is maximal, then its
layer-lengths are in increasing order.  
\end{proof}
\end{lem}

\begin{lem}\label{lem:sizes}
If $|\sigma| \geq 2+\beta$, then we can assume that $\sigma$ begins with an antilayer of size at least $2$ and ends with a 
layer of size at least $\beta$.
\begin{proof}
By Lemmas \ref{lem:layerform} and \ref{lem:sorting}, if $L_k$ is not at least as big as $\beta$, or if $A_1$ is
not at least as big as 2, then there can be {\em{no}} occurrences of $\tau_{2,\beta}$
in $\sigma$.  Assuming $\sigma$ is $\tau_{2,\beta}$-maximal, this cannot be the case.
\end{proof}
\end{lem}

\begin{lem}\label{lem:AnoMax}
We may assume that $|L_k| \geq |A_1|$.
\begin{proof}  
Suppose that $|A_1| > |L_k|$,
and suppose we remove the last element $x$ of $A_1$ and place it at the beginning of $L_k$.
We will show that there is a (strict) net increase in the number of occurrences of $\tau_{2,\beta}$,
so that $\sigma$ cannot have been $\tau_{2,\beta}$-maximal.

First, we count the occurrences which are lost when we remove $x$:
\begin{enumerate}
\item[(1)] if x is a first element, then the second element is in some $L_i$ and the
layer is in some $L_j$, $j>i$: the number of these is
	$$\sum_{1 \leq i < j\leq k} |L_i|\binom{|L_j|}{\beta}.$$
\item[(2)] if x is a second element, then the first element must have been an earlier element
of $A_1$, and the layer could be in any $L_i$: these are enumerated by
	$$(|A_1|-1)\sum_{i=1}^k \binom{|L_i|}{\beta}.$$
Of course $x$ is never in a $\beta$-layer because $x \in A_1$.
\end{enumerate}

Next, we count the occurrences which are gained by putting $x$ into $L_k$.
In all new occurrences, $x$ must be in a new $\beta$-layer, the other elements of which can be any
$\beta-1$ elements from $L_k$.  The first two elements must come from one of three places:
\begin{enumerate}
\item[(1)] both from $A_1$: there are
	$$\binom{|A_1|-1}{2} \binom{|L_k|}{\beta-1}$$
of these.
\item[(2)] one from $A_1$, and one from another layer: there are
	$$(|A_1|-1)\binom{|L_k|}{\beta-1}\sum_{i=1}^{k-1} |L_i|$$
of these.
\item[(3)] from two different layers $L_i$ and $L_j$: there are
	$$\binom{|L_k|}{\beta-1} \sum_{1 \leq i < j \leq k-1} |L_i||L_j|$$
of these.
\end{enumerate}
Hence the loss is
	$$\left( |A_1|-1 \right)\sum_{i=1}^k \binom{|L_i|}{\beta} + \sum_{1 \leq i < j \leq k} |L_i| \binom{|L_j|}{\beta}$$
and the gain is
\begin{equation}\label{eq:thegain}
	\left( |A_1|-1 \right)\sum_{i=1}^{k-1} |L_i|\binom{|L_k|}{\beta-1} 
		+ \binom{|A_1|-1}{2} \binom{|L_k|}{\beta-1} 
		+ \sum_{1 \leq i < j \leq k-1} |L_i||L_j| \binom{|L_k|}{\beta-1}.\;
\end{equation}
Naturally, to provide a contradiction to $\tau_{2,\beta}$-maximality of $\sigma$,
we will be showing that the gain must exceed the loss.  
It will later expedite our computation to consider separately the cases of occurrences of $\tau_{2,\beta}$ which have a $\beta$-layer
in $L_k$:  we can write the loss as the 4-term sum
\begin{equation}\label{eq:bigloss}
	\left(|A_1|-1\right) \sum_{i=1}^{k-1} \binom{|L_i|}{\beta} + \left(|A_1|-1\right) \binom{|L_k|}{\beta}
		+ \sum_{1 \leq i < j \leq k-1} |L_i| \binom{|L_j|}{\beta} + \sum_{i=1}^{k-1} |L_i| \binom{|L_k|}{\beta}.\;
\end{equation}
Now, notice that for each $j \in [1..k]$ we have
\begin{equation}\label{eq:bcrelation}
	\binom{|L_j|}{\beta} = \frac{|L_j|}{\beta}\binom{|L_j|-1}{\beta-1} \leq \frac{|L_k|}{\beta} \binom{|L_k|}{\beta-1}.\;
\end{equation}
because of our assumption that the layers are ordered by increasing size.
Now, we have
\begin{equation*}
\begin{split}
	\left(|A_1|-1\right) \sum_{i=1}^{k-1} \binom{|L_i|}{\beta} &\leq \left( |A_1|-1 \right) \sum_{i=1}^{k-1} \frac{|L_j|}{\beta} \binom{|L_k|}{\beta-1} \\
		&= \frac{1}{\beta} \left[ \left( |A_1|-1 \right) \sum_{i=1}^{k-1} |L_i| \binom{|L_k|}{\beta-1} \right]\;
\end{split}
\end{equation*}
and
\begin{equation*}
\begin{split}
	\sum_{i=1}^{k-1} |L_i|\binom{|L_k|}{\beta} &\leq \sum_{i=1}^{k-1} |L_i| \frac{|L_k|}{\beta} \binom{|L_k|}{\beta-1} \qquad \text{(by \eqref{eq:bcrelation})}\\
		&=\frac{1}{\beta} \left[ |L_k| \sum_{i=1}^{k-1} |L_i| \binom{|L_k|}{\beta-1} \right] \\
		&\leq \frac{1}{\beta} \left[ \left( |A_1|-1 \right)\sum_{i=1}^{k-1} |L_i| \binom{|L_k|}{\beta-1} \right] \;
\end{split}
\end{equation*}
because we are assuming also that $|L_k| < |A_1|$.  Hence, we are able to bound the sum of the first and last terms of \eqref{eq:bigloss}
strictly below the first term of \eqref{eq:thegain} because of our assumption that $\beta \geq 3$.
That the third term of \eqref{eq:bigloss} is bounded below the third term
of \eqref{eq:thegain} is clear from our relation \eqref{eq:bcrelation} on binomial coefficients.
Finally, we may bound the second term of \eqref{eq:bigloss} 
by raw computation: By Lemma \ref{lem:sizes}, we may assume that $|L_k| \geq \beta \geq 3$, 
so that $|A_1| \geq 4$ and thus $|A_1|-1 \leq \frac{3}{2} \left( |A_1|-2 \right)$.
It follows now that
\begin{equation*}
\begin{split}
	\left( |A_1|-1\right) \binom{|L_k|}{\beta} &\leq \left( |A_1|-1\right) \frac{|L_k|}{\beta} \binom{|L_k|}{\beta-1} \\
		&\leq \left( |A_1|-1 \right) \frac{|L_k|}{3} \binom{|L_k|}{\beta-1} \\
		&\leq \left( |A_1|-1 \right) \frac{(|A_1|-1)}{3} \binom{|L_k|}{\beta-1} \\
		&\leq \left( |A_1|-1 \right) \frac{(|A_1|-2)}{2} \binom{|L_k|}{\beta-1} \\
		&= \binom{|A_1|-1}{2} \binom{|L_k|}{\beta-1}.\;
\end{split}
\end{equation*}
We have thus bounded the loss strictly below the gain, so the change must have
resulted in an increase in the number of occurrences of $\tau$.  In particular,
to assume that $|A_1| > |L_k|$ would be in contradiction to the
$\tau_{2,\beta}$-maximality of $\sigma$.
We may assume, then, that $|L_k|$ is at least as big as $|A_1|$.
\end{proof}
\end{lem}

\begin{figure}[b]
\begin{center}

\begin{picture}(120,120)(20,20) 

  \definecolor{dw04blue}{cmyk}{1,.5,0,.1}
  \definecolor{dw04purple}{cmyk}{.3,.9,0,.2}
  \definecolor{dw04orange}{cmyk}{0,.4,1,0}
  \definecolor{dw04black}{cmyk}{0,0,0,1}

	\color{dw04black}
	\put(0,120){\line(1,0){120}}
	\put(0,120){\line(0,-1){120}}
	\put(0,0){\line(1,0){120}}
	\put(120,120){\line(0,-1){120}}

	\color{dw04blue}
	\put(28,118){\line(0,-1){116}}

	\color{dw04blue}
	\put(59,118){\line(0,-1){116}}

	\color{dw04blue}
	\put(2,26){\line(1,0){115}}

	\color{dw04blue}
	\put(2,49){\line(1,0){115}}

	\color{dw04orange}
	\put(3,3){\line(6,5){24}}

	\color{dw04orange}
	\put(32,46){\line(3,-2){25}}

	\color{dw04purple}
	\put(11,10){$A_1$}

	\color{dw04purple}
	\put(39,33){$L_1$}

	\color{dw04purple}
	\put(86,78){$\sigma^*$}

\end{picture} 

\end{center}
\caption{Less-specific structure of the permutation $\sigma$}
\label{fig:alstig}
\end{figure}
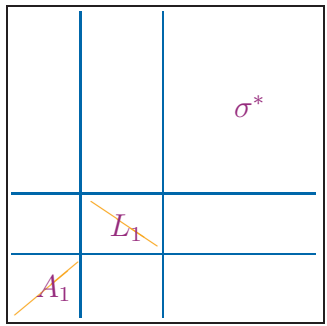

Finally, we may begin piecing together the information we have gathered about $\sigma$ to gain an
important result:

\begin{thm}\label{thm:main}
For $\tau \in S_m$, let 
\begin{equation}\label{eq:layerNoInc}
	g_k(\tau,n) = \max \big\{ g(\tau,\sigma) :: \sigma \in S_n \text{ is of the form } A_1 L_1 \dots L_k \text{ and } |L_k| \geq |A_1|\big\},\;
\end{equation}
that is, maximize only over the permutations $\sigma$ having $k$ nontrivial layers
which satisfy Lemma \ref{lem:AnoMax}.  Then,
for each $k \geq 2$, we have
	$$g_k(\tau_{2,\beta},n) \leq g_1(\tau_{2,\beta},n).$$
\begin{proof}
It will suffice to show that $g_k(\tau_{2,\beta},n) \leq g_{k-1}(\tau_{2,\beta},n)$
for $k \geq 2$.
Let $\sigma = A_1 L_1 \dots L_k \in S_n$ be a permutation for which the expression
\eqref{eq:layerNoInc} is maximized, and assume $k \geq 2$.    
In this context, we will not need the full strength of
Lemma \ref{lem:layerform}; simply write $\sigma = A_1 L_1 \stigma$, as in Figure
\ref{fig:alstig}.  We will show that replacing $A_1 L_1$ with a single antilayer
of size $|A_1| + |L_1|$ does not decrease the number of occurrences of $\tau$.

Let $\Lambda$ denote the number of
occurrences of $\tau_{2,\beta}$ that we lose when we remove $A_1$ and $L_1$, and
let $\Gamma$ be the number of occurrences we gain when we add in the antilayer.
Now, an occurrence of $\tau_{2,\beta}$ is lost whenever it has at least one element in $A_1 \union L_1$.
In this case, either it can have its first two elements in $A_1$ and a layer of size $\beta$ in
$L_1$, its first two elements in $A_1$ and a layer of size $\beta$ in $\sigma^*$, its first two
elements from $A_1$ and $L_1$, respectively, and a layer of size $\beta$ in $\sigma^*$, or
a single element in $A_1 \union L_1$ and the rest of the occurrence (i.e. an occurrence
of $\tau_{1,\beta}$) in $\sigma^*$.  Accordingly, we have 
	$$\Lambda = \binom{|A_1|}{2} \binom{|L_1|}{\beta} + g(\tau_{0,\beta},\sigma^*)\left[ \binom{|A_1|}{2} + |A_1||L_1| \right] + g(\tau_{1,\beta},\sigma^*) \left[ |A_1| + |L_1| \right].$$
When we place an antilayer of size $|A_1|+|L_1|$ at the beginning of $\sigma^*$, we create two
kinds of occurrences.  First, any two elements in the new antilayer, together with any
layer of size $\beta$ in $\sigma^*$, will create an occurrence of $\tau_{2,\beta}$.  Also,
any element of the new antilayer in concert with an occurrence of $\tau_{1,\beta}$ in $\sigma^*$
will create a new occurrence of $\tau_{2,\beta}$.  It follows that
	$$\Gamma = g(\tau_{0,\beta},\sigma^*) \binom{|A_1|+|L_1|}{2} + g(\tau_{1,\beta},\sigma^*) \left[ |A_1|+|L_1| \right],$$
so that the net gain is
\begin{equation}\label{eq:2netgain}
	\Gamma - \Lambda = g(\tau_{0,\beta},\sigma^*) \left[ \binom{|A_1|+|L_1|}{2} - \binom{|A_1|}{2} - |A_1||L_1| \right] - \binom{|A_1|}{2} \binom{|L_1|}{\beta}\;
\end{equation}
(the terms involving $g(\tau_{1,\beta},\sigma^*)$ cancel out).
Now, in the first term of \eqref{eq:2netgain}, the expression $\binom{|A_1|+|L_1|}{2}$
counts the number of ways to choose 2 elements from $A_1 \union L_1$, and we subtract off
the number $|A_1||L_1|$ of ways of choosing one from each set and also the number $\binom{|A_1|}{2}$
of ways to choose both from $A_1$.  If we don't choose both elements from $A_1$ and we don't
choose one from each set, we must choose 2 elements from $|L_1|$, so the first term
of \eqref{eq:2netgain} becomes
	$$g(\tau_{0,\beta},\stigma)\cdot \binom{|L_1|}{2}.$$
Hence, we have
\begin{equation}\label{eq:sim2ng}
	\Gamma - \Lambda = g(\tau_{0,\beta},\stigma)\binom{|L_1|}{2} - \binom{|A_1|}{2}\binom{|L_1|}{\beta}.\;
\end{equation}
We must show, then, that the expression \eqref{eq:sim2ng} is {\em{always}} nonnegative whenever $\stigma$
is nonempty.

Suppose first that $|L_1| \geq |A_1|$.  Then, we have
\begin{equation*}
\begin{split}
	\Gamma - \Lambda &= \binom{|L_1|}{2} g(\tau_{0,\beta},\stigma) - \binom{|A_1|}{2}\binom{|L_1|}{\beta} \\
		&\geq \binom{|A_1|}{2} g(\tau_{0,\beta},\stigma) - \binom{|A_1|}{2}\binom{|L_1|}{\beta} \\
		&= \binom{|A_1|}{2} \left[ g(\tau_{0,\beta},\stigma) - \binom{|L_1|}{\beta} \right].\;
\end{split}
\end{equation*}
Since $\binom{|A_1|}{2}$ is always positive, it is sufficient to prove that $g(\tau_{0,\beta},\stigma) \geq \binom{|L_1|}{\beta}$;
however, this inequality holds whenever $\stigma$ has at least one layer as large as $L_1$, which follows
from Lemma \ref{lem:sorting}.

In the case that $|L_1| < |A_1|$, we know that $\binom{|A_1|}{2} \binom{|L_1|}{\beta} \leq \binom{|L_1|}{2} \binom{|A_1|}{\beta}$
by our technical Lemma (\ref{lem:tech}),
so we have
\begin{equation*}
\begin{split}
	\Gamma - \Lambda &= \binom{|L_1|}{2} g(\tau_{0,\beta},\stigma) - \binom{|A_1|}{2}\binom{|L_1|}{\beta} \\
		&\geq \binom{|L_1|}{2} g(\tau_{0,\beta},\stigma) - \binom{|L_1|}{2}\binom{|A_1|}{\beta} \\
		&= \binom{|L_1|}{2} \left[ g(\tau_{0,\beta},\stigma) - \binom{|A_1|}{\beta} \right],\;
\end{split}
\end{equation*}
which is nonnegative whenever $\stigma$ has a layer at least as large as $|A_1|$.
However, this fact is guaranteed by our adherance to permutations satisfying
Lemma \ref{lem:AnoMax}.  

We have thus reduced the number of nontrivial layers by 1 without decreasing the
number of occurrences of $\tau_{2,\beta}$, so that if a $\tau_{2,\beta}$-maximal permutation $\sigma$
has $k$ nontrivial layers, it still can have no more occurrences of $\tau_{2,\beta}$ than
a $\tau_{2,\beta}$-maximal permutation having $k-1$ nontrivial layers; the result follows.
\end{proof}
\end{thm}
 
\noindent The main results
of this section now follow as corollaries:

\begin{cor}\label{thm:2beta}
There is a pattern $\sigma \in S_n$ which maximizes the number of
occurrences of $\tau_{2,\beta}$ and which consists of a single antilayer 
followed by a single layer.
Hence, the maximum number of occurrences of $\tau_{2,\beta}$ in a permutation in $S_n$ is
\begin{equation}\label{eq:main}
	g(\tau,n) = \max_{x_n \in [0..n]} \binom{x_n}{2} \binom{n-x_n}{\beta}.\;
\end{equation}
\begin{proof}
By Lemmas \ref{lem:layerform} and \ref{lem:AnoMax}, we have 
	$$g(\tau_{2,\beta},n) = \max_{k \leq n} g_k(\tau_{2,\beta},n).$$
From here, it follows from Theorem \ref{thm:main} that $g(\tau_{2,\beta},n) = g_1(\tau_{2,\beta},n)$,
which is clearly \eqref{eq:main}.
\end{proof}
\end{cor}

\begin{cor}
The packing density of $\tau_{2,\beta}$ is 
	$$\rho(\tau_{2,\beta}) = \binom{\beta+2}{2} \left( \frac{2}{\beta + 2} \right)^2 \left( \frac{\beta}{\beta+2} \right)^\beta.$$
\begin{proof}
Notice that we could easily rewrite the expression \eqref{eq:main} as 
\begin{equation}\label{eq:01trans}
	g(\tau,n) = \max_{\xi_n \in [0,1]} \binom{\ceil{\xi_n n}}{2} \binom{\floor{(1-\xi_n)n}}{\beta}.\;
\end{equation}
The work of Price (\cite{price}, Theorem 3.1) shows that a sequence $(\xi_n)$ which maximizes
the expression \eqref{eq:01trans} approaches a constant $\xi$, that is, there is an asymptotically
best ratio of the sizes of the layer and the antilayer of $\sigma$.  
For large enough $n$, we approximate $\xi_n$ by the constant value $\xi$.
It follows that 
\begin{equation}\label{eq:xi}
\begin{split}
	\rho(\tau_{2,\beta}) &= \lim_n \frac{g(\tau_{2,\beta},n)}{\binom{n}{\beta+2}} \\
		&= \lim_n \frac{\binom{ \ceil{\xi n}}{2} \binom{\floor{(1-\xi)n}}{\beta}}{\binom{n}{\beta+2}} \\
		&= \lim_n \frac{ \frac{(\xi n)^2}{2!} \frac{((1-\xi)n)^\beta}{\beta!} }{ \frac{n^{\beta+2}}{(\beta+2)!} } \\
		&= \binom{\beta+2}{2} (1-\xi)^\beta \xi^2.\;
\end{split}
\end{equation}
We may now maximize the expression \eqref{eq:xi} via elementary calculus.  Noting first that
\eqref{eq:xi} evaluates to 0 when $\xi$ is $0$ or $1$, we may maximize the expression 
over $\xi \in [0,1]$ by simply setting its derivative (with respect to $\xi$) to 0: setting
\begin{equation*}
	\frac{\partial}{\partial \xi} \left[ \xi^2 (1-\xi)^\beta \right] = \xi\left[ 2(1-\xi)^\beta - \beta \xi (1-\xi)^{\beta-1} \right] =0,\;
\end{equation*}
we have
	$$2(1-\xi) = \beta \xi$$
since the solutions $\xi = 0$ and $\xi=1$ are clearly unfavorable, so that $\xi = \frac{2}{\beta+2}$.
The result follows.
\end{proof}
\end{cor}


\end{document}